\documentclass[a4paper,11pt]{article}
\usepackage{amssymb}
\usepackage{amsfonts}
\usepackage{amsmath}
\usepackage{a4wide}
\usepackage{jeep}
\usepackage[square]{natbib}
\usepackage{graphicx,slashbox}

\providecommand{\U}[1]{\protect\rule{.1in}{.1in}}
\newtheorem{theorem}{Theorem}

\newtheorem{corollary}[theorem]{Corollary}

\newtheorem{definition}[theorem]{Definition}
\newtheorem{example}[theorem]{Example}

\newtheorem{lemma}[theorem]{Lemma}

\newtheorem{proposition}[theorem]{Proposition}
\newtheorem{remark}[theorem]{Remark}

\newenvironment{proof}[1][Proof]{\noindent \textbf{#1} }{\  \rule{0.5em}{0.5em}}
\begin{document}

\title{\textbf{Free Infinite Divisibility of Free Multiplicative Mixtures of
the Wigner Distribution}}
\author{Victor P\'{e}rez-Abreu\thanks{%
Part of this work was done while the first author was visiting Keio
University during his frequent visits there. He acknowledges the support and
hospitality of the Mathematics Department of this university. } \\
Department of Probability and Statistics, CIMAT \\
Apdo. Postal 402, Guanajuato Gto. 36000, Mexico \\
pabreu@cimat.mx \and Noriyoshi Sakuma\thanks{%
Part of this work was done while the second author was visiting CIMAT.
He sincerely appreciates the support and hospitality of
CIMAT. He is supported by the Japan Society for the Promotion of Science.} \\
Department of Mathematics, Keio University, \\
3-14-1, Hiyoshi, Yokohama 223-8522, Japan.\\
noriyosi@math.keio.ac.jp}
\date{\today}
\maketitle

\begin{abstract}
Let $I^{\ast }$ and $I^{\boxplus }$ be the classes of all classical
infinitely divisible distributions and free infinitely divisible
distributions, respectively, and let $\Lambda $ be the Bercovici-Pata
bijection between $I^{\ast }$ and $I^{\boxplus }.$ The class type $W$ of
symmetric distributions in $I^{\boxplus }$ that can be represented as free
multiplicative convolutions of the Wigner distribution is studied. A
characterization of this class under the condition that the mixing
distribution is 2-divisible with respect to free multiplicative convolution
is given. A correspondence between symmetric distributions in $I^{\boxplus }$
and the free counterpart under $\Lambda $ of the positive distributions in $%
I^{\ast }$ is established. It is shown that the class type $W$ does not
include all symmetric distributions in $I^{\boxplus }$ and that it does not
coincide with the image under $\Lambda $ of the mixtures of the Gaussian
distribution in $I^{\ast }$. Similar results for free multiplicative
convolutions with the symmetric arcsine measure are obtained. Several
well-known and new concrete examples are presented.

\medskip


\end{abstract}

\newpage
\tableofcontents
\section{Introduction}

Let $\mathcal{P}$ denote the set of all Borel probability measures on $%
\mathbb{R}$ and let $\mathcal{P}_{+}$ and $\mathcal{P}_{s}$ be the sets of
all Borel probability measures with support in $\mathbb{R}_{+}=[0,\infty)$
and of all symmetric Borel probability measures (\textit{i.e.} $%
\mu(B)=\mu(-B)$ for all Borel set $B$ on $\mathbb{R}$), respectively.

The free additive convolution $\mu _{1}\boxplus \mu _{2}$ is a binary
operation from $\mathcal{P}\times \mathcal{P}$ to $\mathcal{P}$ that
describes the spectral distribution of the sum $X+Y$ of two freely
independent non-commutative random variables $X$ and $Y$ with spectral
distributions $\mu _{1}$ and $\mu _{2}$, respectively, see \cite{BeVo93}, 
\cite{HiPe00}, \cite{NiSp06}, \cite{VoDyNi92}.

Free infinite divisibility of probability measures with respect to the free
additive convolution $\boxplus $ has received increasing interest during the
last years, see for example \cite{BeNi09}, \cite{BePa99}, \cite{BeVo93}, 
\cite{ChGo08}, \cite{HiMl09}, \cite{Mlpr09}, \cite{Sk09} and references
therein. A key role in free infinite divisibility, similar to the Gaussian
distribution in classical probability, is played by the Wigner or semicircle
distribution $\mathrm{w}$ on$\ (-2,2)$ defined as 
\begin{equation}
\mathrm{w}(\mathrm{d}x)=\frac{1}{2\pi }\sqrt{4-x^{2}}{\Large 1}_{(-2,2)}%
\mathrm{d}x.  \label{FreeBm}
\end{equation}

On the other hand, if $\mu _{1}$ and $\mu _{2}$ are in $\mathcal{P}_{+},$
the free multiplicative convolution $\mu _{1}\boxtimes \mu _{2}$ is the
spectral distribution of $X^{1/2}YX^{1/2}$, where $X$ and $Y$ are freely
independent positive non-commutative random variables with spectral
distributions $\mu _{1}$ and $\mu _{2}$, respectively; see \cite{BeVo93}.
Recently, free multiplicative convolutions for $\mu _{1}$ in $\mathcal{P}_{+}
$ and $\mu _{2}$ in $\mathcal{P}_{s}$ where considered in \cite{APA09}; see
also \cite{RaSp07} when $\mu _{1},\mu _{2}$ have bounded support. It was
also shown in \cite{APA09} that any free symmetric $\alpha $-stable law $\mu
_{\alpha }$, $0<\alpha <2$, can be written as $\mu _{\alpha }=\sigma _{\beta
}\boxtimes \mathrm{w}$, where $\sigma _{\beta }$ is a free positive $\beta $%
-stable law with $\beta =2\alpha /(2+\alpha )$.\ 

The main purpose of the present paper is to study free multiplicative
convolutions $\lambda \boxtimes \mathrm{w}$ of $\lambda \in \mathcal{P}_{+}$
with the Wigner distribution $\mathrm{w}$. In particular, we are interested
in a characterization of the class of free type $W$ distributions consisting
of the measures $\lambda \boxtimes \mathrm{w}$ that are infinitely divisible
with respect to the free additive convolution $\boxplus $. \ 

The analogue problem in classical probability is the study of variance
mixtures of the Gaussian distribution $VZ$, where $V>0$ and $Z$ are
independent (classical) random variables, with $Z$ normally distributed.
Their product is distributed according to the classical multiplicative
convolution of the laws of $V$ and $Z$. It is well-known that $VZ$ is
infinitely divisible in the classical sense if $V^{2}$ is infinitely
divisible and the law of $VZ$ is called of \textit{type }$G$, see\textit{\ }%
\cite{Ke71}, \cite{Ro91}.\textit{\ }However, there are variance mixtures of
the Gaussian distribution which are infinitely divisible but $V^{2}$ is not 
\cite{Ke71}. The relevance of classical type $G$ distributions is that they
are the distributions (at fixed time) of L\'{e}vy processes obtained by
subordination of the classical one-dimensional Brownian motion.

A well-known analytic characterization of type $G$ distributions is as
follows. Recall that a probability measure $\mu $ on $\mathbb{R}$ belongs to
the class $I^{\ast }$ of all infinitely divisible distributions if and only
if the logarithm of its Fourier transform, the so called the \textit{%
classical cumulant transform}, has the L\'{e}vy-Khintchine representation 
\begin{equation}
\mathcal{C}_{\mu }^{\ast }(t)=ib_{\mu }t-\frac{1}{2}{}a_{\mu }t^{2}+\int_{%
\mathbb{R}}\left( \mathrm{e}^{itx}\!-1-itx{\Large 1}_{\left\vert
x\right\vert \leq 1}\right) \nu _{\mu }(\mathrm{d}x),\quad t\in \mathbb{R},
\label{LevyKhintRep}
\end{equation}%
where $b_{\mu }\in \mathbb{R},a_{\mu }\geq 0$, and $\nu _{\mu }$ is a
measure on $\mathbb{R}$ (called the \textit{L\'{e}vy measure}) satisfying $%
\nu _{\mu }(\{0\})=0$ and $\int_{\mathbb{R}}\min (1,x^{2})\nu _{\mu }(%
\mathrm{d}x)<\infty .$ The triplet $(a_{\mu },\nu _{\mu },b_{\mu })$ is
called the \textit{classical generating triplet} of $\mu \in I^{\ast }.$ We
refer to the book by Sato \cite{Sa99} for the study of classical infinitely
divisible distributions on $\mathbb{R}^{d}$. A distribution $\sigma \in
I^{\ast }$ is concentrated on $\mathbb{R}_{+}$ if and only if it admits the 
\textit{regular L\'{e}vy-Khintchine} representation 
\begin{equation}
\mathcal{C}_{\sigma }^{\ast }(t)=ib_{\sigma }t\!+\int_{\mathbb{R}_{+}}\left( 
\mathrm{e}^{itx}\!-1\right) \nu _{\sigma }(\mathrm{d}x),\quad t\in \mathbb{R}%
,  \label{ddslkr}
\end{equation}%
where $b_{\sigma }\geq 0$, $\nu _{\sigma }(-\infty ,0)=0$ and $\int_{\mathbb{%
R}_{+}}(1\wedge x)\nu ($\textrm{$d$}$x)<\infty $. Hence, if $\mu $ is the
distribution of a type $G$ random variable $VZ$, where $V^{2}$ has an
infinitely divisible distribution $\sigma $, it holds that 
\begin{equation}
\mathcal{C}_{\mu }^{\ast }(t)\text{ }=\mathcal{C}_{\sigma }^{\ast
}(it^{2}/2),\quad t\in \mathbb{R}.  \label{ClaCumtypeG}
\end{equation}

The organization and the main results of the paper are as follows. In
Section 2 we present preliminaries and notation on free additive $\boxplus $
and free multiplicative $\boxtimes $ convolutions, and the class \textrm{$I$}%
$^{\boxplus }$ of all $\boxplus $-infinitely divisible distributions. In
particular we present the class of positive regular distributions in \textrm{%
$I$}$^{\boxplus }$. Section 3 introduces the concept of $\boxtimes $--2
divisible distributions in $\mathcal{P}_{+}$ and several examples and
counter-examples of $\boxtimes $--2 divisible distributions are presented.
It is shown that the free Poisson distributions $\mathrm{m}_{c},c\geq 1,$
are $\boxtimes $--2 divisible, but $\mathrm{m}_{c}$ is not for $c<1$ small
enough.

In Section 4 we consider symmetric distribution in \textrm{$I$}$^{\boxplus }$
and their corresponding positive regular distribution in \textrm{$I$}$%
^{\boxplus }$ for which an interesting relation between the corresponding
free cumulant transforms similar to (\ref{ClaCumtypeG}) is proved. This is
surprisingly contrary to the classical case where (\ref{ClaCumtypeG}) does
not hold for all symmetric distributions in $I^{\ast }$. The class of free
type $G$ distributions considered in \cite{ABNPA09} is studied in the
framework of this new relation.

In Section 5 we elaborate about free multiplicative mixtures with the Wigner
distribution. We characterize the class of free type $W$ distributions $%
\overline{\sigma }\boxtimes \mathrm{w,}$ $\overline{\sigma }\in $ $\mathcal{P%
}_{+},$ as those symmetric distributions in \textrm{$I$}$^{\boxplus }$ such
that $\sigma =\overline{\sigma }\boxtimes \overline{\sigma }$ \ is positive
regular in \textrm{$I$}$^{\boxplus }$ and $\boxtimes $--2 divisible. It is
also shown that the distribution in $\mathcal{P}_{+}$ induced by the the
transformation $x\rightarrow x^{2}$ under $\overline{\sigma }\boxtimes 
\mathrm{w}$, $\overline{\sigma }\in $ $\mathcal{P}_{+}$, is always $\boxplus 
$-infinitely divisible and moreover it is a free compound Poisson
distribution. It is shown that there are free type $W$ distributions that
are not free type $G$ distributions and that the Wigner measure is the free
multiplicative convolution of the symmetric arcsine distribution but that
the converse does not hold. Moreover, we prove that the class of free type $%
W $ distributions is not the class of all symmetric free infinitely
divisible distributions.

Finally, in Section 6 we consider symmetric distributions in \textrm{$I$}$%
^{\boxplus }$ which are free multiplicative convolutions with the symmetric
arcsine distribution. We show that this class contains the free type $W$
distributions and the free type $G$ distributions but does not coincide with
the full class of symmetric $\boxplus $-infinitely divisible distributions.

\section{Preliminaries on Free Convolutions and Notation}

In this section we collect some preliminary results and examples on free
infinitely divisible distributions that are used in the remaining of this
work.

\subsection{Free additive convolution $\boxplus$ and infinite divisibility}

Let $\mathbb{C}^{+}$ and $\mathbb{C}^{-}$ be the sets of all complex numbers
satisfying $\mathrm{Im}(z)>0$ and $\mathrm{Im}(z)<0$, respectively. First,
for any probability measure $\mu $ on $\mathbb{R}$, \textbf{the Cauchy
transform} $G_{\mu }:\mathbb{C}^{+}\rightarrow \mathbb{C}^{-}$ is defined by 
\begin{equation*}
G_{\mu }(z)=\int_{\mathbb{R}}\frac{1}{z-x}\mu (\mathrm{d}x),\quad z\in 
\mathbb{C}^{+}.
\end{equation*}%
The reciprocal of the Cauchy transform $F_{\mu }:\mathbb{C}^{+}\rightarrow 
\mathbb{C}^{+}$ of $\mu \in \mathcal{P}$ is defined as $F_{\mu }(z)=\frac{1}{%
G_{\mu }(z)}$. It was shown in \cite{BeVo93} that the right inverse function 
$F_{\mu }^{-1}(z)$ of $F_{\mu }(z)$ (i.e. $F_{\mu }(F_{\mu }^{-1}(z))=z$) is
defined on a region $\Gamma _{\eta ,M}:=\left\{ z\in \mathbb{C};\left\vert 
\mathrm{Re}(z)\right\vert <\eta \mathrm{Im}(z),\quad \mathrm{Im}%
(z)>M\right\} $. This allows us to define \textbf{free cumulant transform}
(or $R$-transform) of a probability measure $\mu $ on $\mathbb{R}$ as $%
\mathcal{C}_{\mu }^{\boxplus }(z)=zF_{\mu }^{-1}(z^{-1})-1$, $z^{-1}\in
\Gamma _{\eta ,M}.$

From the analytic point of view, the \textbf{free additive convolution} of
two probability measures $\mu_{1},\mu_{2}$ on $\mathbb{R}$ is defined as the
probability measure $\mu_{1}\boxplus\mu_{2}$ on $\mathbb{R}$ such that 
\begin{equation}
\mathcal{C}_{\mu_{1}\boxplus\mu_{2}}^{\boxplus}\mathcal{(}z)=\mathcal{C}%
_{\mu_{1}}^{\boxplus}\mathcal{(}z)+\mathcal{C}_{\mu_{2}}^{\boxplus}(z)\quad
z^{-1}\in\Gamma_{\eta,M}  \label{DefAdtCon}
\end{equation}
for $z$ in the common domain where $\mathcal{C}_{\mu_{1}}^{\boxplus}$ and $%
\mathcal{C}_{\mu_{2}}^{\boxplus}$ are defined.

A probability measure $\mu$ on $\mathbb{R}$ is \textbf{free infinitely
divisible} (in short $\boxplus$--ID) if for any $n\in\mathbb{N}$ there
exists a probability measure $\mu_{1/n}$ on $\mathbb{R}$ such that $\mu=\mu
_{1/n}\boxplus\cdot\cdot\cdot\boxplus\mu_{1/n}$ ($n$ times). We denote the
class of all $\boxplus$--ID distributions by $I^{\boxplus}$. As in the
classical case, there is a \textbf{free L\'{e}vy-Khintchine formula} due to
Bercovici and Voiculescu \cite{BeWa08}. In terms of (\ref{DefAdtCon}), $%
\mu\in I^{\boxplus}$ if and only if 
\begin{equation}
\mathcal{C}_{\mu}^{\boxplus}\mathcal{(}z)=b_{\mu}z+a_{\mu}z^{2}+\int_{%
\mathbb{R}}\left( \frac{1}{1-zx}-1-zx1_{\left[ -1,1\right] }\left( x\right)
\right) \nu_{\mu}\left( \mathrm{d}x\right) ,\quad z\in \mathbb{C}^{-\text{ }%
},  \label{FreeLKR}
\end{equation}
where $b_{\mu}\in\mathbb{R},a_{\mu}\geq0$ and $\nu_{\mu}$, the L\'{e}vy
measure, is such that $\nu_{\mu}(\{0\})=0$ and $\int_{\mathbb{R}%
}(1\wedge|x|^{2})\nu_{\mu}(\mathrm{d}x)<\infty$. The triplet $%
(a_{\mu},\nu_{\mu},b_{\mu})$ is unique. When we consider $\boxplus$--ID
distributions, the free cumulant transform can be defined in the lower half
plane $\mathbb{C}^-$. A probability measure $\mu$ is symmetric if and only
if the L\'{e}vy measure $\nu_{\mu}$ is symmetric, $b_{\mu}=0$ and 
\begin{equation}
\mathcal{C}_{\mu}^{\boxplus}\mathcal{(}z)=a_{\mu}z^{2}+\int_{\mathbb{R}%
}\left( \frac{1}{1-zx}-1\right) \nu_{\mu}\left( \mathrm{d}x\right) ,\quad
z\in\mathbb{C}^{-\text{ }}.  \label{FreeLKRSym}
\end{equation}

The \textbf{Bercovici-Pata bijection} $\Lambda:\mathit{I}^{\ast}\rightarrow%
\mathit{I}^{\boxplus}$ between classical and free infinitely divisible
distributions was introduced in \cite{BePa99}. It is such that if $\mu\in%
\mathit{I}^{\ast}$ has classical triplet $(a_{\mu},\nu_{\mu},b_{\mu})$ then $%
\Lambda(\mu)\in\mathit{I}^{\boxplus}$ has free triplet $(a_{\mu},\nu_{%
\mu},b_{\mu})$.

In this work we consider free multiplicative convolutions with the following
key examples. The \textbf{Wigner or semicircle distribution }$\mathrm{w}%
_{b,a}$, with parameters $-\infty <b<\infty ,a>0$ is defined as 
\begin{equation*}
\mathrm{w}_{b,a}(\mathrm{d}x)=\frac{1}{2\pi a}\sqrt{4a-(x-b)^{2}}1_{[b-\sqrt{%
4a},b+\sqrt{4a}]}(x)\mathrm{d}x
\end{equation*}%
and has the free cumulant transform $\mathcal{C}_{\mathrm{w}%
_{b,a}}^{\boxplus }(z)=az^{2}+bz.$ The parameters $b$ and $a$ are the mean
and the variance of this distribution, respectively. It is such that $%
\mathrm{w}_{b,a}=\Lambda (\gamma _{b,a})$, where $\gamma _{b,a}$ is the
classical Gaussian distribution with mean $b$ and variance $a.$ For this
reason $\mathrm{w}_{b,a}$ is also called the \textbf{free Gaussian
distribution.} Especially, we simply write $\mathrm{w}=\mathrm{w}_{0,1}$,
which is corresponding to the standard Gaussian distribution.

Another important example of a free infinitely divisible distribution is the 
\textbf{Marchenko-Pastur distribution} $\mathrm{m_{c}}$ with parameter $c>0$%
, given by 
\begin{equation}
\mathrm{m}_{c}(\mathrm{d}x)=\max(0,(1-c))\delta_{0}(\mathrm{d}x)+\frac{1}{%
2\pi x}\sqrt{4c-(x-1-c)^{2}}{\Large 1}_{[(1-\sqrt{c})^{2},(1+\sqrt{c}%
)^{2}]}(x)\mathrm{d}x.  \label{MaPadis}
\end{equation}
It holds that $\mathrm{m}_{c}=\Lambda(\mathrm{p}_{c})$ where \textrm{p}$_{c}$
is the classical Poisson distribution of mean $c>0$ and it has free triplet $%
(c,c\delta_{1},0)$, where $\delta_{1}$ is the (probability) measure
concentrated at one. For this reason $\mathrm{m}_{c}$ is also called the 
\textbf{free Poisson distribution.} In the case $c=1$, we simply write $%
\mathrm{m}=\mathrm{m}_{1}.$ 

We now consider the case of free infinitely divisible distributions with
nonnegative support, for which we have a situation different than for the
classical case (\ref{ddslkr}), where the drift is nonnegative, the L\'{e}vy
measure is concentrated on $\mathbb{R}_{+}$ and there is not Gaussian part.
In the free case we consider two situations. First, following a similar
terminology as in \cite{PR07}, we propose to call a distribution $\sigma \in
I_{+}^{\boxplus }$ \textbf{free regular} or simply regular, if its L\'{e}vy
Khintchine representation is given by 
\begin{equation}
\mathcal{C}_{\sigma }^{\boxplus }\mathcal{(}z)=b_{\sigma }z+\int_{\mathbb{R}%
_{+}}\left( \frac{1}{1-zx}-1\right) \nu _{\sigma }\left( \mathrm{d}x\right)
,\quad z\in \mathbb{C}^{-\text{ }},  \label{FreeRLKR}
\end{equation}%
where $b_{\sigma }\geq 0,$ $\nu _{\sigma }((-\infty ,0])=0$ and $%
\int_{0}^{\infty }\min (1,x)\nu _{\sigma }\left( \mathrm{d}x\right) <\infty
. $ Not all nonnegative free infinitely divisible distributions are regular.
For example, the standard Wigner distribution translated by $2$, that is $%
\mathrm{w}_{2,1}(\mathrm{d}x)=\Lambda (\gamma _{2,1})(\mathrm{d}x)=\frac{1}{%
2\pi }\sqrt{4-(x-2)^{2}}1_{[0,4]}(x)\mathrm{d}x$ has support on $(0,4),$ its
free triplet is $(1,0,2)$ and hence there is a non-zero Wigner part $a_{%
\mathrm{w}_{2,1}}=1$. Let $I_{r+}^{\boxplus }$ be the class of all regular
distribution in $I^{\boxplus }.$ It holds that $\Lambda (\mathit{I}%
_{+}^{\ast })=$\textrm{$I$}$_{r+}^{\boxplus }$, but $\mathit{I}%
_{r+}^{\boxplus }$ is not equal to $\{\mu \in \mathit{I}^{\boxplus }\,\,|\,\,%
\text{the support of}\,\,\mu \,\,\text{is in}\,\,\mathbb{R}_{+}\}$. The free
Poisson distribution for all $c>0$ is free regular infinitely divisible. 

When a distribution $\mu $ has all moments, the \textbf{free cumulants} $%
\kappa _{n}$ are the coefficients of the formal expansion 
\begin{equation*}
\mathcal{C}_{\mu }^{\boxplus }(z)=\sum_{n=1}^{\infty }\kappa _{n}z^{n}.
\end{equation*}%
The following is a easy necessary criteria for free infinite divisibility
which is used repeatedly in this work. It is based on the first four
cumulants of a distribution. It is the analogue to a criterion in classical
infinite divisibility in \cite[pp 181]{SV03}. It was pointed out to us by O.
Arizmendi. Its proof follows from the fact that the classical cumulant
sequence $(c_{n})$ of $\mu \in \mathit{I}^{\ast }$ (i.e. the coefficients $%
(c_{n})$ of the series expansion 
\begin{equation*}
\mathcal{C}_{\mu }^{\ast }(t)=\sum_{n=1}^{\infty }\frac{c_{n}}{n!}%
t^{n},\quad t\in \mathbb{R},
\end{equation*}%
of the classical cumulant transform) coincide with the free cumulant
sequence $(\kappa _{n})$ of $\Lambda (\mu )$. For other criteria using free
cumulants see the recent paper \cite{Wl09}. The book \cite{NiSp06}
elaborates in general conditions related to the moment problem.

\begin{lemma}
\label{kurtosis} \label{lemAr} Let $\mu$ be an $\boxplus$--ID measure having
four finite moments. Then $\kappa_{2}\kappa_{4}\geq\kappa_{3}^{2}$.
\end{lemma}

An example of a distribution that is not free infinitely divisible but plays
an important role in this work is the \textbf{symmetric arcsine distribution}
with parameter $s$%
\begin{equation}
\mathrm{a}_{s}(\mathrm{d}x)=\frac{1}{\pi}\frac{1}{\sqrt{s-x^{2}}}{\LARGE 1}%
_{(-\sqrt{s},\sqrt{s})}(x)\mathrm{d}x.  \label{arcsinmea}
\end{equation}
If the parameter $s=1$, we use the notation $\mathrm{a}$ for $\mathrm{a}_{1}$%
. It is interesting that the distribution arises as the additive convolution 
$\mathrm{a}=\mathrm{d\boxplus d}$, where $\mathrm{d}$ is the symmetric
Bernoulli (atomic) measure $\mathrm{d}=\frac{1}{2}\left( \delta_{
-1}+\delta_{1}\right) $.


\subsection{Free multiplicative convolution $\boxtimes$ and the S-transform}

For a probability measure $\mu$ on $\mathbb{R}$, the $p$--th push--forward
measure of $\mu^{(p)}$ of $\mu$ is defined as 
\begin{equation*}
\mu^{(p)}(B)=\int_{\mathbb{R}}1_{B}(\left\vert x\right\vert ^{p})\mu(\mathrm{%
d}x),\quad B\in\mathcal{B}((0,\infty)).
\end{equation*}
It is trivial to see that $\mathrm{w}^{(2)}=\mathrm{m}$ and $\mathrm{a}%
^{(2)}=\mathrm{a}^{+}$, where $\mathrm{a}^{+}$ is is the positive arcsine
law on $(0,1)$ given by \eqref{postarcsin} below.

To study the free \lq\lq product\textquotedblright\ $\boxtimes$ of
probability measures, it is useful to consider another analytic tool called
the \textbf{S-transform, }which is defined as follows. For $\mu\in\mathcal{P}%
_{+}$ define 
\begin{equation}
\Psi_{\mu}(z)=\int_{\mathbb{R}}\frac{zx}{1-zx}\mu(\mathrm{d}x)=z^{-1}G_{\mu
}(z^{-1})-1,\quad z\in\mathbb{C}\backslash\mathbb{R}.  \label{S-Cauchy}
\end{equation}

It was proved in \cite{BeVo93} that for probability measures with support on 
$\mathbb{R}_{+}$ and such that $\mu(\{0\})<1,$ the function $\Psi_{\mu}(z)$
has a unique inverse $\chi_{\mu}(z)$ in the left-half plane $i\mathbb{C}^{+}$
and $\Psi_{\mu}(i\mathbb{C}^{+})$ is a region contained in the circle with
diameter $(\mu(\{0\})-1,0)$. In this case the $S$-transform of $\mu$ is
defined as $S_{\mu}(z)=\chi_{\mu}(z)\frac{1+z}{z}$ It satisfies $z=\mathcal{C%
}_{\mu}^{\boxplus}\left( zS_{\mu}(z\right) )$ for sufficiently small $%
z\in\Psi_{\mu}(i\mathbb{C}^{+})$.

Following \cite{BeVo93}, the\textbf{\ free multiplicative convolution} of a
probability measure $\mu_{1}$, $\mu_2$ supported on $\mathbb{R}_{+}$ is
defined as the positive probability measure $\mu_{1}\boxtimes\mu_{2}$ on $%
\mathbb{R}$ such that 
\begin{equation}
S_{\mu_{1}\boxtimes\mu_{2}}\mathcal{(}z)=S_{\mu_{1}}\mathcal{(}z)S_{\mu_{2}}%
\mathcal{(}z)  \label{DefMulConPos}
\end{equation}
for $z$ in a common region of $\Psi_{\mu_{1}}(i\mathbb{C}^{+})\cup\Psi
_{\mu_{2}}(i\mathbb{C}^{+})$.

The definition of $S$-transform was extended for symmetric probability
measures $\mu\in\mathcal{P}_{s}$ in \cite{APA09}. Let $H=\left\{ z\in\mathbb{%
C}^{-};\quad\left\vert \mathrm{Re}(z)\right\vert <|\mathrm{Im}(z)|\right\}
,\quad\widetilde{H}=\left\{ z\in \mathbb{C}^{+};\quad\left\vert \mathrm{Re}%
(z)\right\vert <\mathrm{Im}(z)\right\}$. It was proved in \cite{APA09} that
when $\mu\in\mathcal{P}_{s}$ with $\mu(\{0\})<1$, the transform $\Psi_{\mu}$
has a unique inverse on $H$, $\chi_{\mu}:\Psi_{\mu}(H)\rightarrow H$ and a
unique inverse on $\widetilde{H}$, $\widetilde{\chi}_{\mu}:\Psi_{\mu}(%
\widetilde{H})\rightarrow\widetilde{H}.$ In this case there are two $S$%
-transforms for $\mu$ given by 
\begin{equation}  \label{DefStran}
S_{\mu}(z)=\chi_{\mu}(z)\frac{1+z}{z}\text{ and }\widetilde{S}_{\mu }(z)=%
\widetilde{\chi}_{\mu}(z)\frac{1+z}{z}
\end{equation}
\medskip and these are such that%
\begin{equation}
S_{\mu}^{2}(z)=\frac{1+z}{z}S_{\mu^{(2)}}(z)\text{ and }\widetilde{S}_{\mu
}^{2}(z)=\frac{1+z}{z}S_{\mu^{(2)}}(z)  \label{STSym}
\end{equation}
for $z$ in $\Psi_{\mu}(H)$ and $\Psi_{\mu}(\widetilde{H})$, respectively.
Moreover the following result holds.

\begin{lemma}
Assume that $\mu\in\mathcal{P}_{s}\cup\mathcal{P}_+$. For some sufficiently
small $\varepsilon>0$, we have a region $D_{\varepsilon}$ that includes $%
\{-it;0<t<\varepsilon\}$ such that 
\begin{equation}
z=\mathcal{C}_{\mu}^{\boxplus}(zS_{\mu}(z))  \label{critical}
\end{equation}
for $z\in D_{\varepsilon}$.
\end{lemma}

\begin{proof}
For $\mu \in \mathcal{P}_+$ see \cite{NiSp97}. Let $\mu \in \mathcal{P}_s$.
We take some $\alpha>0$ and $\beta>0$ such that there exists $%
F_{\mu}^{-1}(z) $ for $z\in\Gamma_{\alpha,\beta}$. Then we have the inverse
in $\{z\in \mathbb{C}^{-};|z|>\alpha\Im z,\,\,|z|<1/\beta\}.$ From $%
\Psi_{\mu}(z)=\Psi_{\mu^{(2)}}(z^{2})$, and \cite[Proposition 6.1]{BeVo93}, $%
\lim _{it\rightarrow0,t<0}\Psi_{\mu}(it)=0$. Then $\lim_{t\rightarrow0,t<0}%
\tilde{\chi}_{\mu}(t)=0$. If we take sufficiently small $\varepsilon>0$, $z%
\tilde{S}_{\mu}(z)=(z+1)\tilde{\chi}_{\mu}(z)$ maps from sufficiently small
domain that contains $\{-it;0<t<\varepsilon^{\prime }<\varepsilon\}$ to $%
\{z\in\mathbb{C}^{-};|z|>\alpha\Im z,\,\,|z|<1/\beta\}$ for some $%
\varepsilon^{\prime}$ smaller than $\varepsilon$. Then we always have some
region where we have \eqref{critical} from \cite{APA09}.
\end{proof}

Following \cite{APA09}, the\textbf{\ free multiplicative convolution} of a
probability measure $\mu_{1}$ supported on $\mathbb{R}_{+}$ with a symmetric
probability measure $\mu_{2}$ on $\mathbb{R}$ is defined as the symmetric
probability measure $\mu_{1}\boxtimes\mu_{2}$ on $\mathbb{R}$ such that 
\begin{equation}
S_{\mu_{1}\boxtimes\mu_{2}}(z)=S_{\mu_{1}}(z)S_{\mu_{2}}(z).
\label{DefMulCon}
\end{equation}
It was also shown in \cite{APA09} that 
\begin{equation}
\mu_{1}\boxtimes\mu_{2}^{(2)}\boxtimes\mu_{1}=(\mu_{1}\boxtimes%
\mu_{2})^{(2)}.  \label{numunu}
\end{equation}

For the convenience of the reader, we include other examples of the $S$%
-transform of important examples, which appear repeatedly in this paper. For
the Wigner measure $\mathrm{w}_{0,a}$ with zero mean and variance $a$ 
\begin{equation}
S_{\mathrm{w_{0,a}}}(z)=\sqrt{\frac{1}{az}}  \label{WigmeasStrans}
\end{equation}%
For the Marchenko-Pastur measure $\mathrm{m_{c}}$ with parameter $c>0$ 
\begin{equation}
S_{\mathrm{m}_{c}}(z)=\frac{1}{z+c}.  \label{MarchPasturStrans}
\end{equation}

For the arcsine distribution (\ref{arcsinmea}) 
\begin{equation}
S_{\mathrm{a}_{s}}(z)=\sqrt{\frac{z+2}{sz}.}  \label{arcsinStrans}
\end{equation}
For the \textbf{positive arcsine distribution} on $(0,s)$%
\begin{equation}
\mathrm{a}_{s}^{+}(\mathrm{d}x)=\frac{1}{\pi}\frac{1}{\sqrt{x(s-x)}}{\LARGE 1%
}_{(0,s)}(x)\mathrm{d}x  \label{postarcsin}
\end{equation}
its $S$-transform is 
\begin{equation}
S_{\mathrm{a}_{s}^{+}}(z)=\frac{z+2}{s(z+1)}.  \label{postarcsinStrans}
\end{equation}
When $s=1$, we use the notation $\mathrm{a}^{+}$for $\mathrm{a}_{1}^{+}$.


\subsection{$\boxplus$-compound Poisson distributions}

Using the Bercovici-Pata bijection we define \textbf{free compound Poisson
distributions}. For a combinatorial treatment of the free compound Poisson
distribution see the book \cite{NiSp06}.

\begin{definition}
Let $\sigma$ be a probability measure with $\sigma(\{0\})=0$ and let $c$ be
positive number. \newline
a) $\mu$ is compound Poisson distribution $(c,\sigma)$ if its cumulant
transform can be represented as 
\begin{equation*}
\mathcal{C}_{\mu}^{\ast}(t)=c(\hat{\sigma}(t)-1)=c\int_{\mathbb{R}%
}(\exp(itx)-1)\sigma(\mathrm{d}x),\quad t\in\mathbb{R}.
\end{equation*}
b) $\mu$ is free compound Poisson distribution $(c,\sigma)$ on $\mathbb{R}$
if $\Lambda^{-1}(\mu)$ is a classical compound Poisson distribution $%
(c,\sigma)$. In this case 
\begin{equation}
\mathcal{C}_{\mu}^{\boxplus}(z)=c\int_{\mathbb{R}_{+}}\left( \frac{1}{1-zx}%
-1\right) \sigma(\mathrm{d}x)\quad z\in\mathbb{C}^{-}.  \label{freecumCP}
\end{equation}
\end{definition}

We denote by $\mu^{\boxplus c}$ the $c$ times free additive convolution of $%
\mu\in\mathcal{P}$.

\begin{proposition}
\label{FCPResults} a) If a probability measure $\mu$ in \textrm{$I$}$%
_{r+}^{\boxplus}$ or \textrm{$I$}$_{s}^{\boxplus}$ is free compound Poisson
distribution $(c,\sigma)$, then $\mu^{\boxplus1/c}=m\boxtimes\sigma$.

b) If $\mu=\mathrm{m}\boxtimes\sigma$ for some $\sigma$ in $\mathcal{P}_{+}$
or $\mathcal{P}_{s}$ respectively, then $\mu$ is the free compound Poisson
distribution.

c) If $\sigma\in\mathcal{P}_{+}$, then $\mu=\mathrm{m}\boxtimes\sigma\in $%
\textrm{$I$}$_{r+}^{\boxplus}.$
\end{proposition}

\begin{proof}
$a)$ Since $\mu$ is free infinitely divisible, for any $c>0$ the measure $%
\mu^{\boxplus1/c}$ exists. Therefore we have that $\mathcal{C}%
_{\mu^{\boxplus1/c}}^{\boxplus}(z)=\int_{\mathbb{R}_{+}}\left( \frac{1}{1-zx}%
-1\right) \sigma(\mathrm{d}x)$ for $z\in\mathbb{C}^{-}$. Now, using (\ref%
{S-Cauchy}) we have $\int_{\mathbb{R}}(\frac{zx}{1-zx})\sigma(\mathrm{d}%
x)=\int_{\mathbb{R}}(\frac{1}{1-zx}-1)\sigma(\mathrm{d}x)=\Psi_{\sigma}(z)$.
Then 
\begin{equation*}
\mathcal{C}_{\mu^{\boxplus1/c}}^{\boxplus}(zS_{\mu^{\boxplus1/c}}(z))=z=%
\Psi_{\sigma}(zS_{\mu^{\boxplus1/c}}(z)),\quad\text{for sufficiently small }%
z\in\Psi_{\mu^{\boxplus1/c}}(i\mathbb{C}^{+}),
\end{equation*}
and therefore $\chi_{\sigma}(z)=zS_{\mu^{\boxplus1/c}}(z)$. Thus from %
\eqref{DefStran} and \eqref{MarchPasturStrans}, $S_{\mu^{\boxplus1/c}}(z)=S_{%
\mathrm{m}}(z)S_{\sigma}(z)$ and we conclude that $\mu^{\boxplus1/c}=m%
\boxtimes\sigma$. \newline
b) Assume that $\mu=m\boxtimes\nu$. The S-transform of $\mu$ is $S_{\mu}(z)=%
\frac {1}{1+z}S_{\nu}(z)=\frac{1}{z}\chi_{\nu}(z)$. Therefore $%
\Psi_{\nu}(zS_{\mu }(z))=z=\mathcal{C}_{\mu}^{\boxplus}(zS_{\mu}(z))$. We
have $\mathcal{C}_{\mu }^{\boxplus}(z)=\Psi_{\nu}(z)=\int_{\mathbb{R}}(\frac{%
zx}{1-zx})\nu (\mathrm{d}x)=\int_{\mathbb{R}}(\frac{1}{1-zx}-1)\nu(\mathrm{d}%
x)$. We then conclude that $\mu=\sigma\boxtimes\nu$ is free compound Poisson
distribution. \newline
(c) From (b), the L\'{e}vy measure of $\mu$ is $\sigma$. So, $\int_{\mathbb{R%
}_{+}}\min(1,x)\sigma(\mathrm{d}x)<\infty$ and $\sigma$ is concentrated on $%
\mathbb{R}_{+}$. Therefore $\mu$ is \textrm{$I$}$_{r+}^{\boxplus}$.
\end{proof}


\section{$\boxtimes$--2 Divisibility of Probability Measures}

\label{2div}

In this section we consider the concept of $\boxtimes$--2 divisibility that
it is used in Section 5 to characterize free type $W$ distributions.

\begin{definition}
A probability measure $\sigma\in\mathcal{P}_{+}$ is called $\boxtimes $-%
\textbf{2 divisible} if there exists a probability measure $\overline {\sigma%
}\in\mathcal{P}_{+}$ such that $\sigma=\overline{\sigma}\boxtimes \overline{%
\sigma}$.
\end{definition}

Bercovici and Voiculescu \cite{BeVo93} consider the more general concept of
infinite divisibility of probability measures on $\mathbb{R}_{+}$ with
respect to the free multiplicative convolution \textbf{$\boxtimes$}. Namely,
a probability measure $\mu$ on $\mathbb{R}_{+}$ is $\boxtimes$--infinitely
divisible if for any $n\in\mathbb{N}$, there exists $\mu_{n}\in \mathcal{P}%
_+ $ such that 
\begin{equation*}
\mu=\underbrace{\mu_{n}\boxtimes\mu_{n}\boxtimes\cdots\boxtimes\mu_{n}}%
_{n\,\,\text{times}}.
\end{equation*}
Of course, any \textbf{$\boxtimes$}-infinitely divisible distribution is $%
\boxtimes$--2 divisible. For results on indecomposable measures we refer to
the recent paper \cite{BeWa082}.

Examples of free multiplicative infinitely divisible distributions are the
following.

\begin{example}
\label{bxtinf} (1) $\delta_{0}$ is $\boxtimes$--infinitely divisible. 
\newline
(2) The free Poisson distribution $\mathrm{m}_{c}$ is $\boxtimes$%
--infinitely divisible if and only if $c\geq1.$ This is because if 
\begin{equation*}
\Sigma_{\mathrm{m}_{c}}(z)=S_{\mathrm{m}_{c}}\left( \frac{z}{1-z}\right)
=\exp\left( -\log\frac{1-z}{(1-c)z+c}\right) ,
\end{equation*}
the function $v(z):=-\log\frac{1-z}{(1-c)z+c}$ is analytic on $\mathbb{C}%
\backslash\mathbb{R}_{+}$, $v(\overline{z})=\overline{v(z)}$ and $v(\mathbb{C%
}^{+})\subset\mathbb{C}^{-}$ if $c\geq1$, and not analytic on $\mathbb{C}%
\backslash\mathbb{R}_{+}$ if $c<1$. Hence Theorem 6.13. in \cite{BeVo93}
implies that $m_c$ is $\boxtimes$-infinitely divisible if and only if $c\geq
1$ (see also Theorem 1.2. in \cite{BBCC09}). \newline
(3) The positive $\boxplus$--stable laws, with index $0<\alpha<1$, are also $%
\boxtimes$--infinitely divisible. This follows from \cite[Proposition A 4.4]%
{BePa99}.
\end{example}

We are able to obtain a class of examples of $\boxtimes$--2 divisible
distributions, from distributions using Fuss-Catalan numbers $A_{m}(p,r)$
recently constructed in M\l otkowski \cite{Mlpr09}. The Fuss-Catalan numbers
are defined for, $p\in\mathbb{R}$ and $r\in\mathbb{R}$, as follows: $%
A_{0}(p,r)=1$ and 
\begin{equation*}
A_{m}(p,r):=\frac{r}{m!}\prod_{i=1}^{m-1}(mp+r-i)\quad\text{if}\quad m\geq1.
\end{equation*}
It was proved in \cite{Mlpr09} that when $p\geq1$ and $0\leq r\leq p$, the
Fuss-Catalan numbers $\{A_{m}(p,r)\}_{m=1}^{\infty}$ are the moments of a
probability measure concentrated on $\mathbb{R}_{+}$, denoted by $%
\mu_{(p,r)} $. We call $\mu_{(p,r)}$ the M\l otkowski distribution with
parameter $p,r$ and the set of all such measures the M\l otkowski class.
Furthermore the following was shown in \cite{Mlpr09}.

\begin{lemma}
\label{MlotDis}(a) The free cumulant sequence of $\mu_{(p,r)}$ is $%
\{A_{m}(p-r,r)\}_{m=1}^{\infty}$. \newline
(b) If $0\leq2r\leq p$ and $r+1\leq p$ then $\mu_{(p,r)}$ is $\boxplus$%
-infinitely divisible. \newline
(c) $\mu_{(p_{1},r)}\boxtimes\mu_{(1+p_{2},1)}=\mu_{(p_{1}+rp_{2},r)}$ for $%
r\not =0$. \newline
\end{lemma}

With the above lemma we can construct examples of $\boxtimes$--2 divisible
distributions and consider their $\boxplus$-infinitely divisibility.

\begin{example}
\label{ExamPoi2Div} (1) Since all the cumulants of the Marchenko-Pastur
distribution $\mathrm{m}=\mathrm{m}_{1}$ are equal to one, from (a) and (c)
in the above example we have that $\mathrm{m}=\overline{\mathrm{m}}\boxtimes%
\overline{\mathrm{m}}$ is $\boxtimes$--2 divisible with $\overline{\mathrm{m}%
}=\mu_{(3/2,1)}$.

(2) It is easy to see that $\mu_{(3/2,1)}$ does not satisfy the condition in
Lemma \ref{lemAr}. Then $\mu_{(3/2,1)}$ is not $\boxplus$--infinitely
divisible.

(3) From Lemma \ref{lemAr} we have that the measure $\mu_{(5/4,1)}$ is not $%
\boxplus$--infinitely divisible. This is because the first free cumulants of 
$\mu_{(5/4,1)}$ are $\kappa_{2}=1/2$, $\kappa_{3}=3/24$ and $k_{4}=0$ and
they do not satisfy the condition in Lemma \ref{lemAr}. However, the
distribution $\mu_{(3/2,1)}=\mu_{(5/4,1)}\boxtimes\mu_{(5/4,1)}$ is $%
\boxtimes$--2 divisible but not $\boxtimes$--infinitely divisible.
\end{example}

From Example \ref{bxtinf} we have the following result.

\begin{proposition}
\label{2DivMarPas} For $c\geq1$ the free Poisson distribution $\mathrm{m}%
_{c} $ is $\boxtimes$--2 divisible.
\end{proposition}

On the other hand, for $0<c<1$ sufficiently small we can prove that $\mathrm{%
m}_{c}$ is not $\boxtimes $-2 divisible.

\begin{proposition}
\label{No2divMP}For $0<c<1$ sufficiently small the function $S_{m_c}(z)=%
\frac{1}{\sqrt{z+2c}}$ is not the $S$-transform of a probability measure on $%
\mathbb{R}_{+}$.
\end{proposition}

\begin{proof}
Assume there exists a probability measure $\sigma$ with the $S$--transform $%
S_{\sigma}(z)=\frac{1}{z+c}.$ From definition of the $S$--transform (\ref%
{DefStran}) we have 
\begin{equation*}
\Psi_{\sigma}\left( \frac{z}{1+z}S_{\sigma}(z)\right) =\Psi_{\sigma}\left( 
\frac{z}{(1+z)\sqrt{z+c}}\right) =z.
\end{equation*}
Consider the expansion 
\begin{equation*}
\frac{z}{(1+z)\sqrt{z+c}}=\frac{z}{\sqrt{c}}+\frac{(-2c-1)z^{2}}{2c^{3/2}}+%
\frac{\left( 8c^{2}+4c+3\right) z^{3}}{8c^{5/2}}+\frac{\left(
-16c^{3}-8c^{2}-6c-5\right) z^{4}}{16c^{7/2}}+O\left( z^{5}\right) ,
\end{equation*}
the inverse-series satisfies 
\begin{equation*}
\Psi_{\sigma}(z)=\sqrt{c}z+\frac{1}{2}(2c+1)z^{2}+\frac{\left(
8c^{2}+12c+1\right) z^{3}}{8\sqrt{c}}+\left( c^{2}+3c+1\right) z^{4}+O\left(
z^{5}\right) .
\end{equation*}
So the first three moments of $\sigma$ are $m_1(\sigma) = \sqrt{c}, \quad
m_2(\sigma) = \frac{1}{2}(2c+1),\quad m_3(\sigma) = \frac{8c^{2}+12c+1}{8%
\sqrt{c}}$. Then $\det(m_{i+j}(\sigma))_{0\leq i\leq2,0\leq j\leq2}=\frac{%
8c-1}{64t}$. If we take $c<1/8$, $\det(m_{i+j}(\sigma))_{0\leq i\leq2,0\leq
j\leq2}<0$, which is a contradiction to the Stieltjes moment problem.
\end{proof}

From the above examples we have the following summary which is a useful
result for the study of free type $W$ distributions in Section 5.

\begin{proposition}
(1) Let $\overline{\sigma} = \mu_{(3/2,1)}$. Then $\overline{\sigma}\not \in 
\mathit{I}_{r+}^{\boxplus}$ but $\sigma=\overline{\sigma}\boxtimes\overline{%
\sigma}\in\mathit{I}_{r+}^{\boxplus}$. \newline
(2) Let $\overline{\sigma} = \mu_{(5/4,1)}$. Then $\overline{\sigma}\not \in 
\mathit{I}_{r+}^{\boxplus}$ and $\sigma=\overline{\sigma}\boxtimes \overline{%
\sigma}\not \in \mathit{I}_{r+}^{\boxplus}$. \newline
(3) $\{\overline{\sigma}\in\mathcal{P}_{+};\overline{\sigma}\boxtimes 
\overline{\sigma}\in\mathit{I}_{r+}^{\boxplus}\}\not =\mathcal{P}_{+}.$
\end{proposition}

It is an open problem to find a free regular probability measure $\overline{%
\sigma}\in I_{r+}^{\boxplus}$ such that $\sigma=\overline{\sigma }\boxtimes%
\overline{\sigma}$ is not $\boxplus$--infinitely divisible. We have the
following table

\begin{table}[th]
\caption{$\boxtimes$--2 divisibility and $\boxplus$--regular ID}
\begin{center}
\begin{tabular}{|c||c|c|}
\hline
\backslashbox{$\sigma$}{$\overline{\sigma}$} & $\mathit{I}^{\boxplus}_{r+}$
& \makebox{not $\mathit{I}^{\boxplus}_{r+}$} \\ \hline\hline
$\mathit{I}^{\boxplus}_{r+}$ & $\mathrm{m}^{\boxtimes2}=\mathrm{m}\boxtimes%
\mathrm{m}$ and $\boxplus$--stable case & $\mathrm{m} = \mu
_{(3/2,1)}\boxtimes\mu_{(3/2,1)}$ \\ \hline
not $\mathit{I}^{\boxplus}_{r+}$ & No example & $\mu_{(3/2,1)} = \mu
_{(5/4,1)}\boxtimes\mu_{(5/4,1)}$ \\ \hline
\end{tabular}%
\end{center}
\end{table}


\bigskip

\section{Symmetric $\boxplus$--Infinitely Divisible Distributions}

In this section we prove an interesting relation between the free cumulant
transform of a symmetric free infinitely divisible distributions and that of
an associated regular positive free infinitely divisible distribution. As
to our knowledge, there is not such a relation in the classical case. We
consider the particular case of the free type $G$ distribution studied in 
\cite{ABNPA09} as the image of the classical type $G$ distributions under
the Bercovici-Pata bijection $\Lambda.$ We prove that free type $G$
distributions can be represented as mixtures of free multiplicative
convolutions.

Using the push-forward notation of a measure, for a probability measure $\mu$
on $\mathbb{R}_{+}$, let $\mu^{(1/2)+}$ and $\mu^{(1/2)-}$ be the induced
measures by $\mu$ on $(0,\infty)$ and $(-\infty,0)$ under the mappings $%
\sqrt{x}\rightarrow x$ and $-\sqrt{x}\rightarrow x$, respectively. We
observe that $\mu^{(1/2)+}$ and $\mu^{(1/2)-}$ are L\'{e}vy measures if $%
\mu$ is L\'{e}vy measure of a free regular infinitely divisible distribution.

\subsection{A characterization}

\begin{theorem}
\label{main1} $\mu\in$\textrm{$I$}$_{s}^{\boxplus}$ if and only if there is $%
\sigma\in$\textrm{$I$}$_{r+}^{\boxplus}$ such that 
\begin{equation}
\mathcal{C}_{\mu}^{\boxplus}(z)=\mathcal{C}_{\sigma}^{\boxplus}(z^{2}).\quad
z\in\mathbb{C}\backslash\mathbb{R}.  \label{dec1}
\end{equation}
Moreover, the L\'{e}vy measures of $\mu$ and $\sigma$ are related by a symmetrization
\begin{equation}
\nu_{\mu}=\frac{1}{2}\left(
\nu_{\sigma}^{(1/2)+}+\nu_{\sigma}^{(1/2)-}\right)  \label{dec}
\end{equation}
and 
\begin{equation}
\nu_{\sigma}=2\nu_{\mu}^{(2)}.  \label{dec01}
\end{equation}
\end{theorem}

\begin{remark} 
a) Since $\mathcal{C}_{\mathrm{w}}^{\boxplus }(z)=z^{2}$
for the Wigner distribution \textrm{w}, 
from (\ref{dec1}) we have that for any $%
\mu \in $\textrm{$I$}$_{s}^{\boxplus }$, it holds that there is $\sigma \in $%
$I_{r+}^{\boxplus }$ such that $\mathcal{C}_{\mu }^{\boxplus }(z)=%
\mathcal{C}_{\sigma }^{\boxplus }(\mathcal{C}_{\mathrm{w}}^{\boxplus
}(z)),\,z\in \mathbb{C}\backslash \mathbb{R}$. In the classical case, we have
from (\ref{ClaCumtypeG}) that if $\gamma _{1}$ is the standard classical
Gaussian distribution, it holds that  only for for some $\mu ^{\prime }\in 
I_{s}^{\ast },$ $\mathcal{C}_{\mu ^{\prime }}^{\ast }(t)$ $=%
\mathcal{C}_{\sigma ^{\prime }}^{\ast }(i\mathcal{C}_{\gamma _{1}}^{\ast
}(t)),t\in \mathbb{R}$, with $\sigma ^{\prime }\in $\textrm{$I$}$_{+}^{\ast }
$. 

b) However, the classical corresponding of (\ref{dec1}) is that for any $\Lambda^{-1}(\mu) \in I_{s}^{\ast }$ there is $\Lambda^{-1}(\sigma)\in $%
\textrm{$I$}$_{+}^{\ast }$ such that   
\begin{equation*}
\mathcal{C}_{\Lambda^{-1}(\mu)}^{\ast }(t)=\int_{\mathbb{R}_{+}}\left( \cos (t%
\sqrt{x}-1\right) \nu _{\Lambda^{-1}(\sigma)}(\mathrm{d}x),\quad
t\in \mathbb{R}\text{.}
\end{equation*}

c) For illustration of this theorem, consider the standard Wigner
distribution $\mu =\mathrm{w}$ for which $\mathcal{C}_{\mathrm{w}}^{\boxplus
}(z)=z^{2}$. In this case $\sigma =\delta _{1}$ and $\mathcal{C}_{\sigma
}^{\boxplus }(z)=z$. For the standard Cauchy distribution $\mu $ we have $%
\mathcal{C}_{\mu }^{\boxplus }(z)=-iz$ and $\nu _{\mu }(dx)=\frac{1}{x^{2}}%
1_{\mathbb{R}}(x)dx$. In this case $\mathcal{C}_{\sigma }^{\boxplus }(z)=-i%
\sqrt{z}$ is the free cumulant transform of the one-side $1/2$-free stable
distribution with L\'{e}vy measure $\nu _{\sigma }(dx)=\frac{1}{x^{3/2}}1_{%
\mathbb{R}_{+}}(x)dx$.
\end{remark}

\begin{proof}
If we have $\sigma\in I^{\boxplus}_{+}$, $\mathcal{C}^{\boxplus}_{\sigma}(z)$
can be define on $\mathbb{C}\backslash\mathbb{R}_{+}$. Suppose $\mu\in $%
\textrm{$I$}$_{s}^{\boxplus}$ and let $z\in\mathbb{C}\backslash\mathbb{R}$.
Then, 
\begin{align*}
\mathcal{C}_{\mu}^{\boxplus}(z) & =a_{\mu}z^{2}+\int_{\mathbb{R}}\left( 
\frac{1}{1-zx}-1\right) \nu_{\mu}(\mathrm{d}x) \\
& =a_{\mu}z^{2}+\int_{\mathbb{R}}\left( \frac{1}{1-zx}-1-zx1_{|x|\leq
1}(x)\right) \nu_{\mu}(\mathrm{d}x) \\
& =a_{\mu}z^{2}+\int_{\mathbb{R}_{+}}\left( \frac{1}{1-zx}-1-zx1_{|x|\leq
1}(x)\right) \nu_{\mu}(\mathrm{d}x) \\
& +\int_{\mathbb{R}\backslash\mathbb{R}_{+}}\left( \frac{1}{1-zx}%
-1-zx1_{|x|\leq1}(x)\right) \nu_{\mu}(\mathrm{d}x).
\end{align*}
Since the L\'{e}vy measure $\nu_{\mu}$ is symmetric, we have 
\begin{align}
\mathcal{C}_{\mu}^{\boxplus}(z) & =a_{\mu}z^{2}+\int_{\mathbb{R}_{+}}\left( 
\frac{1}{1-zx}-1-zx1_{|x|\leq1}(x)\right) \nu_{\mu}(\mathrm{d}x)  \notag \\
& +\int_{\mathbb{R}_{+}}\left( \frac{1}{1+zx}-1+zx1_{|x|\leq1}(x)\right)
\nu_{\mu}(\mathrm{d}x)  \notag \\
& =a_{\mu}z^{2}+\int_{\mathbb{R}_{+}}\left( \frac{1}{1-zx}+\frac{1}{1+zx}%
-2\right) \nu_{\mu}(\mathrm{d}x)  \notag \\
& =a_{\mu}z^{2}+2\int_{\mathbb{R}_{+}}\left( \frac{1}{1-z^{2}x^{2}}-1\right)
\nu_{\mu}(\mathrm{d}x)  \notag \\
& =a_{\mu}z^{2}+2\int_{\mathbb{R}_{+}}\left( \frac{1}{1-z^{2}x}-1\right)
\nu_{\mu}^{(2)}(\mathrm{d}x).  \label{rel1}
\end{align}
Let $\nu_{\sigma}=$ $2\nu_{\mu}^{(2)}$. Then $\nu_{\sigma}((-\infty,0])=0$
and using again the symmetry of $\nu_{\mu}$ we have 
\begin{align*}
\int_{\mathbb{R}_{+}}\min(1,x)\nu_{\sigma}\left( \mathrm{d}x\right) & =\int_{%
\mathbb{R}_{+}}\min(1,x)2\nu_{\mu}^{(2)}\left( \mathrm{d}x\right) \\
& =2\int_{\mathbb{R}}\min(1,x^{2})\nu_{\mu}\left( \mathrm{d}x\right) <\infty,
\end{align*}
since $\nu_{\mu}$ is a L\'{e}vy measure. Then $\nu_{\sigma}$ is the L\'{e}vy
measure of a free regular infinitely divisible distribution $\sigma$ and
using (\ref{rel1}) and the uniqueness of the L\'{e}vy-Khintchine
representation we have $\mathcal{C}_{\mu}^{\boxplus}(z)=\mathcal{C}%
_{\sigma}^{\boxplus}(z^{2})$. Here $\sigma\in\mathit{I}_{r+}^{\boxplus}$ has
triplet $(0,\nu_{\sigma },c_{\sigma})$ with $c_{\sigma}=a_{\mu}$.

Conversely, if a $\sigma \in \mathit{I}_{r+}^{\boxplus }$ with triplet $%
(0,\nu _{\sigma },c_{\sigma })$, let $\mu $ be the symmetric free infinitely
divisible with triplet $(a_{\sigma },\nu _{\sigma },0)$ where $a_{\mu
}=c_{\mu }$ and $\nu _{\mu }=\frac{1}{2}\left( \nu _{\sigma }^{(1/2)+}+\nu
_{\sigma }^{(1/2)-}\right) $. Then $\nu _{\mu }$ is a symmetric L\'{e}vy
measure and 
\begin{align*}
\mathcal{C}_{\mu }^{\boxplus }(z)& =a_{\mu }z^{2}+\int_{\mathbb{R}}\left( 
\frac{1}{1-zx}-1\right) \nu _{\mu }(\mathrm{d}x) \\
& =a_{\mu }z^{2}+\frac{1}{2}\int_{\mathbb{R}}\left( \frac{1}{1-zx}-1\right)
\left( \nu _{\sigma }^{(1/2)+}(\mathrm{d}x)+\nu _{\sigma }^{(1/2)-}(\mathrm{d%
}x)\right) \\
& =a_{\mu }z^{2}+\frac{1}{2}\int_{\mathbb{R}_{+}}\left( \frac{1}{1-z\sqrt{x}}%
+\frac{1}{1+z\sqrt{x}}-2\right) \nu _{\sigma }(\mathrm{d}x) \\
& =a_{\mu }z^{2}+\int_{\mathbb{R}_{+}}\left( \frac{1}{1-z^{2}x}-1\right) \nu
_{\sigma }(\mathrm{d}x)=\mathcal{C}_{\sigma }^{\boxplus }(z^{2}),
\end{align*}%
which proves the result.
\end{proof}

The relations between the L\'{e}vy measures (\ref{dec}) and (\ref{dec01})
hold in the free and classical infinitely divisible cases. The importance of
the above theorem is in the relation (\ref{dec1}) between the cumulant
transforms.

When $\sigma$ is absolutely continuous, we have the following formula.

\begin{lemma}
\label{deeq} If the L\'{e}vy measure of $\sigma \in $\textrm{$I$}$_{+}^{\ast
}$ has density $h_{\sigma }$, then the $\nu _{\sigma }^{(1/2)+}$ is also
absolutely continuous and has density $h_{\sigma }^{(1/2)+}(x)$ of $\nu
_{\sigma }^{(1/2)+}$ 
\begin{equation*}
h_{\sigma }^{(1/2)+}(x)=2xh_{\sigma }(x^{2})\quad (x>0).
\end{equation*}
\end{lemma}

\subsection{Classical and free type $G$ distributions}

Following \cite{ABNPA09}, we say that a probability distribution $\upsilon$
is in the class of \textbf{free type }$G$ distributions, if there exists a
classical type $G$ distribution $\mu$ such that $\upsilon=\Lambda(\mu)$.
That is, $\Lambda^{-1}(\upsilon)$ is the distribution of $VZ$ where $V$ and $%
Z$ are independent, with $Z$ having the standard Gaussian distribution and $%
V^{2}$ has a classical infinitely divisible distribution in $I_{+}^{\ast}$
with the L\'{e}vy measure $\rho_{V}.$

The following result is a characterization of the L\'{e}vy measure of
classical and free type $G$ distributions in terms of the measure $\sigma$
of Theorem \ref{main1}. We also consider the relation between $\sigma$ and
the distribution of $V^{2}$ which is given in terms of mixtures of free
multiplicative convolutions.

\begin{theorem}
A non Gaussian symmetric probability measure $\mu$ on $\mathbb{R}$ is a type 
$G$ distribution iff the L\'{e}vy measure of $\sigma$ can be represented as 
\begin{equation*}
\nu_{\sigma}(\mathrm{d}x)=\frac{g(x)}{\sqrt{x}}\mathrm{d}x
\end{equation*}
where $g(x)$ is a completely monotone function on $(0,\infty)$.
\end{theorem}

\begin{proof}
If $\mu$ is type $G$ distribution then its L\'{e}vy measure can be
represented $\nu_{\mu}(\mathrm{d}x)=g(x^{2})$ where $g(x)$ is a completely
monotone function on $(0,\infty)$. 
\begin{equation*}
\nu_{\mu}(\mathrm{d}x)=g(x^{2})\mathrm{d}x=\frac{1}{2}\left(
2g(x^{2})1_{(0,\infty)}(x)+2g(x^{2})1_{(-\infty,0)}(x)\right) \mathrm{d}x
\end{equation*}
So, $\nu_{\sigma}^{(1/2)+}=2g(x^{2})1_{(0,\infty)}(x)\mathrm{d}x$ and $%
\nu_{\sigma}^{(1/2)-}=2g(x^{2})1_{(-\infty,0)}(x)\mathrm{d}x$. By Lemma \ref%
{deeq}$,\nu_{\sigma}(\mathrm{d}x)=h_{\sigma}(x)\mathrm{d}x=\frac {g(x)}{%
\sqrt{x}}\mathrm{d}x$. The converse is trivial.
\end{proof}

It is well known that the L\'{e}vy measure of a type $G$ distribution $\mu$
is of the form $\nu_{\mu}(\mathrm{d}x)=v_{\mu}(x)\mathrm{d}x$ where%
\begin{equation}
v_{\mu}(x)=\int_{\mathbb{R}_{+}}\phi(x,s)\rho_{V}(\mathrm{d}s)
\label{LevMeaTypeG}
\end{equation}
for some L\'{e}vy measure $\rho_{V}$ of a distribution in $I_+^{\ast}$ and $%
\phi(x,s)$ is the Gaussian density of mean zero and variance $s$.

We have seen that a free type $G$ distribution $\mu$ is related to two L\'{e}%
vy measures on $\mathbb{R}_{+};$ $\nu_{\sigma}$ given by Theorem \ref{main1}
and $\rho_{V\text{ }}$ as in (\ref{LevMeaTypeG})$.$ The following result
gives the relation between these two L\'{e}vy measures in terms of mixtures
of cumulants of free multiplicative convolutions of the Marchenko-Pastur
distribution with $\gamma_{s}^{(2)},$ the Gamma distribution with shape
parameter $1/2$ and scale parameter $s.$

\begin{theorem}
\label{Main3} Let $\mu $ be a free type $G$ distribution with the L\'{e}vy
measure $\nu _{\mu }$ given by (\ref{LevMeaTypeG}) with the mixing L\'{e}vy
measure $\rho _{V}.$ $\ $Let $\sigma \in $\textrm{$I$}$_{r+}^{\boxplus }$
with the L\'{e}vy measure $\nu _{\sigma }=2\nu _{\mu }^{(2)}$. Then 
\begin{equation}
\mathcal{C}_{\sigma }^{\boxplus }(z)=\int_{\mathbb{R}_{+}}\mathcal{C}_{%
\mathrm{m}\boxtimes \gamma _{s}^{(2)}}^{\boxplus }(z)\rho _{V}(\mathrm{d}s).
\label{TypeGCumMixt2}
\end{equation}%
Moreover, 
\begin{equation}
\nu _{\sigma }(dx)=\int_{\mathbb{R}}\gamma _{s}^{(2)}(dx)\rho _{V}(ds).
\label{Newequat}
\end{equation}
\end{theorem}

\begin{proof}
Using the free regular representation (\ref{FreeRLKR}), the fact that $%
\nu_{\sigma}=2\nu_{\mu}^{(2)}$ and Proposition \ref{FCPResults}(b) with the
free cumulant transform of the compound Poisson distribution $\mathrm{m}%
\boxtimes\gamma_{s}^{(2)}$in (\ref{freecumCP}) we have 
\begin{align*}
\mathcal{C}_{\sigma}^{\boxplus}(z) & =a_{\mu}z+\int_{\mathbb{R}_{+}}\left( 
\frac{1}{1-zx}-1\right) \nu_{\sigma}(\mathrm{d}x) \\
& =a_{\mu}z+2\int_{\mathbb{R}_{+}}\left( \frac{1}{1-zx}-1\right) \nu_{\mu
}^{(2)}(\mathrm{d}x) \\
& =a_{\mu}z+\int_{\mathbb{R}}\left( \frac{1}{1-zx^{2}}-1\right) \nu_{\mu }(%
\mathrm{d}x) \\
& =a_{\mu}z+\int_{\mathbb{R}}\left( \frac{1}{1-zx^{2}}-1\right) \int_{%
\mathbb{R}_{+}}\phi(x,s)\rho_{V}(\mathrm{d}s)\mathrm{d}x\quad(\text{use }%
\eqref{LevMeaTypeG}) \\
& =a_{\mu}z+\int_{\mathbb{R}_{+}}\int_{\mathbb{R}}\left( \frac{1}{1-zx^{2}}%
-1\right) \phi(x,s)\mathrm{d}x\rho_{V}(\mathrm{d}s) \\
& =a_{\mu}z+\int_{\mathbb{R}_{+}}\int_{\mathbb{R}_{+}}\left( \frac{1}{1-zx}%
-1\right) \gamma_{s}^{(2)}(\mathrm{d}x)\rho_{V}(\mathrm{d}s) \\
& =a_{\mu}z+\int_{\mathbb{R}_{+}}\mathcal{C}_{\mathrm{m}\boxtimes\gamma
_{s}^{(2)}}^{\boxplus}(z)\rho_{V}(\mathrm{d}s).
\end{align*}
The second statement follows also from the above calculations.
\end{proof}

If in the above theorem we consider the variance mixture of the Normal
distribution with the classical Poisson distribution, we obtain $\Lambda (%
\mathcal{L}(VZ))=\mathrm{m}\boxtimes\gamma_{1}$. 

\section{Free Type \textit{W} Distributions}

In this section we consider multiplicative convolutions of the Wigner
measure with probability measures on $\mathbb{R}_{+}$ and their free
additive infinite divisibility. We introduce the class of \textbf{free type }%
$W$\textbf{\ distributions} and show the role played by $\boxtimes$-2
divisible distributions.

\subsection{Multiplicative mixtures of the Wigner distribution}

\begin{definition}
For a probability measure $\lambda $ in $\mathcal{P}_{+}$, a symmetric
probability measure $\mu =\lambda \boxtimes \mathrm{w}$ is called a free
multiplicative mixture of the Wigner distribution. When $\mu $ is free
infinitely divisible, it is called a \textbf{free type }$W$\textbf{\
distribution}.
\end{definition}

From \cite[Theorem 12]{APA09} we have that any free symmetric $\alpha $%
-stable law is a free type $W$ distribution. However, not all all symmetric
distribution can be represented as free multiplicative mixture of the Wigner
distribution. This is the case of the arcsine symmetric distribution $%
\mathrm{a}$ on $(-1,1)$, as shown by the following result.

\begin{proposition}
Let $\mathrm{a}$ be arcsine distribution on $(-1,1)$. There does not exist $%
\lambda\in\mathcal{P}_{+}$ such that $\mathrm{a}=\lambda\boxtimes\mathrm{w}$.
\end{proposition}

\begin{proof}
If there is such a $\lambda\in\mathcal{P}_{+}$, from (\ref{DefMulCon}) we
have $S_{\mathrm{a}}(z)=S_{\lambda}(z)S_{\mathrm{w}}(z)$. So, using (\ref%
{WigmeasStrans}) and (\ref{arcsinStrans}) we obtain $S_{\lambda}(z)=\sqrt{z+2%
}$. But $S_{\lambda}^{\prime}(t)=\frac{1}{2\sqrt{t+2}}$ is positive for all $%
t\in(-1,0)$. Then, from Proposition 6.8 in \cite{BeVo93}, $S_{\lambda}(z)=%
\sqrt{z+2}$ cannot be the $S$-transform of a probability measure on $\mathbb{%
R}_{+}$.
\end{proof}

A remarkable property of an arbitrary free multiplicative mixture of the
Wigner measure is the fact that\ "its square" is always free infinitely
divisible and moreover a free compound Poisson distribution.

\begin{proposition}
\label{Wigmultsq} \label{Rmrakble2} Let $\overline{\sigma}\in\mathcal{P}_{+}$%
, $\sigma =\overline{\sigma}\boxtimes\overline{\sigma}$, \textrm{w} be
Wigner measure and $\mu=\overline{\sigma}\boxtimes\mathrm{w}$. Then $%
\mu^{(2)}=\sigma \boxtimes\mathrm{m}$ is free infinitely divisible.
Moreover, $\mu^{(2)}$ is a free compound distribution in \textrm{$I$}$%
_{r+}^{\boxplus}.$
\end{proposition}

\begin{remark}
In the classical case we have a similar result than in the above case.
Namely, any square of a variance mixture of Gaussian (whether $\ast$%
-infinitely divisible or not) is always $\ast$-infinitely divisible. To see
this, let $X=VZ$ for a positive random variable $V$ independent of $Z$ with
the standard Gaussian distribution. Recall the well known fact that $Z$ has
the distribution of $E^{1/2}A$, where $E$ and $A$ are independent random
variables, $E$ with exponential distribution and $A$ with a symmetric
arcsine distribution on $(-1,1)$. Then $X^{2}$ has the distribution of $%
EV^{2}A^{2}$ which is always $\ast$-infinitely divisible, since any mixture
of the exponential distribution is $\ast$-infinitely divisible. Moreover, $%
X^{2}$ is in the Bondesson class of distributions characterized by $\ast$%
-infinitely divisible distributions with completely monotone L\'{e}vy
density, see for example \cite{BMS06}.
\end{remark}

\begin{proof}
Since $\sigma=\overline{\sigma}\boxtimes\overline{\sigma}$ and $\mu =%
\overline{\sigma}\boxtimes\mathrm{w}$, using (\ref{DefMulConPos}), (\ref%
{DefMulCon}) and the fact that $S_{\mathrm{w}}(z)=\frac{1}{\sqrt{z}}$ we
have $S_{\sigma}(z)=S_{\overline{\sigma}}(z)^{2}$ and $S_{\mu}(z)=\frac {1}{%
\sqrt{z}}S_{\overline{\sigma}}(z)$. On the other hand, from (\ref{STSym}) $%
S_{\mu}^{2}(z)=\frac{z+1}{z}S_{\mu^{(2)}}(z)=\frac{1}{z}S_{\sigma}(z)$.
Hence, $S_{\mu^{(2)}}(z)=\frac{1}{z+1}S_{\sigma}(z)$ and since for the
Marchenko-Pastur distribution $\mathrm{m}$ $S_{\mathrm{m}}(z)=\frac{1}{z+1}$%
, again using (\ref{DefMulConPos}) we obtain $\mu^{(2)}=\sigma\boxtimes 
\mathrm{m}$.
\end{proof}

In particular, if we consider the square of symmetric free stable law, the
free compound Poisson distribution with one-side stable law appears.

\begin{corollary}
Let $\upsilon_{\alpha}$ be a symmetric free $\alpha$-stable law, $0<\alpha<2$%
. Then $\upsilon_{\alpha}^{(2)}=(\sigma_{\beta}\boxtimes\sigma_{\beta}
)\boxtimes\mathrm{m}$, where $\sigma_{\beta}$ is a free positive $\beta $%
-stable law with $\beta=2\alpha/(2+\alpha)$.
\end{corollary}

\begin{proof}
The result follows from the last proposition, relation (\ref{numunu}) and
Proposition \ref{FCPResults} (b).
\end{proof}

\subsection{Type $W$ and $\boxtimes$-2 divisibility}

We now present a characterization of type $W$ distributions where the
concept of 2--$\boxtimes$ divisibility appears in connection with free
regular infinitely distributions on $\mathbb{R}_{+}$. Specifically, a free
multiplicative mixture of the Wigner distribution $\mu$ is free infinitely
divisible if and only if the mixing measure $\sigma$ is free regular and 2--$%
\boxtimes$ divisible. Moreover, $\sigma$ is the distribution appearing in
Theorem (\ref{main1}).

\begin{theorem}
\label{main2} \label{main0} Let $\overline{\sigma}\in\mathcal{P}_{+}$ and 
\textrm{w} be the Wigner measure. Then $\sigma=\overline{\sigma}\boxtimes 
\overline{\sigma}\in$\textrm{$I$}$_{r+}^{\boxplus}$ if and only if $\mu=%
\overline{\sigma}\boxtimes\mathrm{w}\in$\textrm{$I$}$_{s}^{\boxplus}$ in
which case 
\begin{equation*}
\mathcal{C}_{\mu}^{\boxplus}(z)=\mathcal{C}_{\sigma}^{\boxplus}(z^{2}),\quad
z\in\mathbb{C}\backslash\mathbb{R}.
\end{equation*}
\end{theorem}

\begin{proof}
Assume first that $\overline{\sigma}\boxtimes\overline{\sigma}\in \mathit{I}%
_{r+}^{\boxplus}$. From theorem \ref{main1} there exists a symmetric $%
\boxplus$--infinitely divisible distribution $\mu$ such that $\mathcal{C}%
_{\mu}^{\boxplus}(z)=\mathcal{C}_{\overline{\sigma}\boxtimes\overline{\sigma}%
}^{\boxplus}(z^{2})$. Then 
\begin{align*}
\mathcal{C}_{\overline{\sigma}\boxtimes\overline{\sigma}}^{\boxplus }(zS_{%
\overline{\sigma}}^{2}(z)) & =\mathcal{C}_{\mu}^{\boxplus}(\sqrt {z}S_{%
\overline{\sigma}}(z)) \\
& =\mathcal{C}_{\mu}^{\boxplus}(zS_{\overline{\sigma}\boxtimes\mathrm{w}}(z))
\\
& =\mathcal{C}_{\mu}^{\boxplus}(zS_{\mu}(z)).
\end{align*}
Therefore, $S_{\mu}(z)=S_{\overline{\sigma}\boxtimes\mathrm{w}}(z)$ and by
the uniqueness of the multiplicative convolution we obtain $\mu=\overline{%
\sigma }\boxtimes\mathrm{w}$.\newline

Next, for a probability measure $\overline{\sigma}$ on $\mathbb{R}_{+}$ let $%
\mu=\overline{\sigma}\boxtimes\mathrm{w}$, that is $S_{\mu}(z)=S_{\overline {%
\sigma}\boxtimes\mathrm{w}}(z)=\frac{1}{\sqrt{z}}S_{\overline{\sigma}}(z)$.
Again, using theorem \ref{main1} there exists a regular $\boxplus$%
--infinitely divisible distribution $\sigma$ (on $\mathbb{R}_{+})$ such that 
$\mathcal{C}_{\mu}^{\boxplus}(z)=\mathcal{C}_{\sigma}^{\boxplus}(z^{2})$.
Then 
\begin{align*}
\mathcal{C}_{\mu}^{\boxplus}(zS_{\mu}(z)) & =\mathcal{C}_{\mu}^{\boxplus }(%
\sqrt{z}S_{\overline{\sigma}}(z)) \\
& =\mathcal{C}_{\mu}^{\boxplus}(zS_{\overline{\sigma}}^{2}(z)) \\
& =\mathcal{C}_{\sigma}^{\boxplus}(zS_{\overline{\sigma}\boxtimes \overline{%
\sigma}}(z)).
\end{align*}
Then, from (\ref{critical}) we have that%
\begin{equation*}
z=\mathcal{C}_{\mu}^{\boxplus}(zS_{\mu}(z))=\mathcal{C}_{\sigma}^{\boxplus
}(zS_{\overline{\sigma}\boxtimes\overline{\sigma}}(z))=\mathcal{C}_{\sigma
}^{\boxplus}(zS_{\sigma}(z))
\end{equation*}
and therefore $S_{\sigma}(z)=S_{\overline{\sigma}\boxtimes\overline{\sigma}%
}(z)$ and hence $\sigma=\overline{\sigma}\boxtimes\overline{\sigma}\in%
\mathit{I}_{r+}^{\boxplus}$.
\end{proof}

\smallskip From Example \ref{ExamPoi2Div}(3) we have that if $\overline {%
\sigma}=\mu_{(5/4,1)}$, then $\sigma=\overline{\sigma}\boxtimes \overline{%
\sigma}=\mu_{(3/2,1)}$ is $\boxtimes$-$2$ divisible but it is not $\boxplus$%
--infinitely divisible. Then, from the above theorem we have that $%
\mu_{(5/4,1)}\boxtimes\mathrm{w}$ is a free multiplicative mixture of the
Wigner distribution but not a type $W$ distribution. So we have the
following result.

\begin{proposition}
The class of all type $W$ distributions is a proper subset of the class of
all free multiplicative mixtures of the Wigner distribution.
\end{proposition}

\subsection{Relation to free type $G$ distributions}

We now give an example of a type $W$ distribution which is not a free type $%
G $ distribution and therefore type $W$ distribution is not the image of the
class of classical type $G$ distributions under the Bercovici-Pata bijection 
$\Lambda.$

Let $\mathrm{b}_{s}$ be the symmetric Beta distribution $(1/2,3/2)$ on $(-2%
\sqrt{s},2\sqrt{s})$ given by 
\begin{equation*}
\mu_{s,\alpha,\beta}(\mathrm{d}x)=\frac{1}{2B(\alpha,\beta)\sqrt{s}}%
|x|^{\alpha-1}(s-|x|)^{\beta-1}1_{(-2\sqrt{s},2\sqrt{s})}(x)\mathrm{d}x.
\end{equation*}
for $\alpha=1/2,\beta=3/2.$ It was shown in \cite{ABNPA09} that $\mathrm{b}%
_{s}$ is $\boxplus$-infinitely divisible distribution which is not a free
type $G$ distribution. Moreover $\mathrm{b}_{s}=\mathrm{a}_{s}\boxtimes%
\mathrm{m}$. Then $S_{\mathrm{b}_{1}}(z)=S_{\mathrm{a}}(z)S_{\mathrm{m}}(z)$
and using (\ref{MarchPasturStrans}) and (\ref{arcsinStrans}) we obtain that 
\begin{equation}
S_{\mathrm{b}_{1}}(z)=\frac{1}{\sqrt{z}}\frac{\sqrt{z+2}}{z+1}.  \label{aux1}
\end{equation}
With this we can prove that $\mathrm{b}_{1}$ is a type $W$ distribution.

\begin{lemma}
\label{SBisTW} The symmetric Beta distribution $\mathrm{b}_{1}$ on $(-2,2)$
is a type $W$ distribution. Moreover $\mathrm{b}_{1}=\mathrm{w}\boxtimes 
\mathrm{a^{+}}\boxtimes\mathrm{\overline{m_{2}}}$ where $\mathrm{a^{+}}$ is
the positive arcsine distribution on (0,1) and $\mathrm{m}_{2}$ is the free
Poisson distribution of parameter $c=2.$
\end{lemma}

\begin{proof}
Using (\ref{aux1}) we have 
\begin{equation*}
S_{\mathrm{b}_{1}}(z)=\frac{1}{\sqrt{z}}\frac{\sqrt{z+2}}{z+1}=\frac{1}{%
\sqrt{z}}\frac{z+2}{z+1}\frac{1}{\sqrt{z+2}}.
\end{equation*}
The result follows, since from (\ref{MarchPasturStrans}) and Proposition \ref%
{2DivMarPas} $S_{\mathrm{\overline{m_{2}}}}(z)=1/\sqrt{z+2},$ and from (\ref%
{postarcsinStrans}) we obtain $S_{\mathrm{a^{+}}}(z)=(z+2)/(z+1).$
\end{proof}

Finally, we show that type $W$ distributions is a proper subclass of the
class \textrm{$I$}$_{s}^{\boxplus}$ of all symmetric free infinitely
divisible distributions$.$ The example below is constructed from the family
of the free Poisson distributions.

\begin{example}
\label{exNW} Let $\widetilde{\mathrm{m}}_{c}$ be dual of $\mathrm{m}_{c}$,
that is $\widetilde{\mathrm{m}}_{c}(B)=\mathrm{m}_{c}(-B)$ for any Borel set 
$B$. Then for $c>0$ 
\begin{equation*}
\widetilde{\mathrm{m}}_{c}(\mathrm{d}x)=\max (0,(1-c))\delta _{0}(\mathrm{d}%
x)+\frac{1}{2\pi (-x)}\sqrt{4c-((-x)-1-c)^{2}}{\Large 1}_{[-(1+\sqrt{c}%
)^{2},-(1-\sqrt{c})^{2}]}(x)\mathrm{d}x.
\end{equation*}%
Let $\mu _{c}=\mathrm{m}_{c}\boxplus \widetilde{\mathrm{m}}_{c}$. Then $\mu
_{c}$ is a symmetric $\boxplus $--infinitely divisible distribution with the
free cumulant transform 
\begin{equation*}
\mathcal{C}_{\mu _{c}}^{\boxplus }(z)=\frac{2cz^{2}}{1-z^{2}},
\end{equation*}%
and the $S$--transform 
\begin{equation}
S_{\mu _{c}}(z)=\frac{1}{\sqrt{z}}\frac{1}{\sqrt{z+2c}}.  \label{aux2}
\end{equation}%
However, for small $c$, $\mu _{c}$ is not a free multiplicative convolution
of the Wigner measure. Indeed, it was shown in Proposition \ref{No2divMP}
that for $c$ small enough (for example $c<1/16)$, the function $1/\sqrt{z+2c}
$ is not the $S$-transform of a probability measure on $\mathbb{R}_{+}.$
Then, since $1/\sqrt{z}$ is the $S$-transform of the Wigner measure, from (%
\ref{aux2}) and (\ref{DefMulCon}) we have that $\mu _{c}$ is not the $S$%
-transform of a multiplicative convolution of the Wigner measure. However,
for $c\geq 1/2,$ $\mathrm{m}_{c}\boxplus \widetilde{\mathrm{m}_{c}}=%
\overline{\mathrm{m}}_{2c}\boxtimes \mathrm{w}.$
\end{example}

We summarize some results of this section as follows.

\begin{theorem}
(1) The class of all type $W$ distributions is a proper subset of \textrm{$I$%
}$_{s}^{\boxplus}$.\newline
(2) The intersection of the class of all type $W$ distributions and the
class of all free type $G$ distributions is not empty.\newline
(3) The class of all type $W$ distributions does not coincide with the class
of all free type $G$ distributions.
\end{theorem}


\section{$\boxplus$--ID of Free Multiplicative Convolutions with the Arcsine
Measure}

We now define a new subclass of $\boxplus$--infinitely divisible
distributions, the \textbf{type} $\mathbf{AS}$ distributions. We say that a
distribution $\mu\in\mathit{I}_{s}^{\boxplus}$, belongs to the class type $%
AS $, if it is the free multiplicative convolution of the arcsine measure.
That is, there exists a distribution $\lambda\in\mathcal{P}_{+}$ such that $%
\mu=\lambda\boxtimes\mathrm{a}$.

It was already seen that $\mathrm{w=\overline{m_{2}}}\boxtimes\mathrm{a}$.
Then the class type $AS$ contains the class of type $W$ distributions. We
characterize the class of type $AS$ distributions in a similar way as the
class type $W$.

\begin{theorem}
A symmetric distribution $\mu=\lambda\boxtimes\mathrm{a,}$ $\lambda \in%
\mathcal{P}_{+}$, is a type $AS$ distribution if and only if there exists $%
\sigma\in\mathit{I}_{r+}^{\boxplus}$ such that $\lambda\boxtimes \lambda=%
\mathrm{m}_{2}\boxtimes\sigma$.
\end{theorem}

\begin{proof}
Assume $\mu=\lambda\boxtimes\mathrm{a}$ is the type $AS$ distribution. Then $%
S_{\mu}(z)=S_{\lambda}(z)S_{\mathrm{a}}(z)$ and from Theorem \ref{main1},
there exists $\sigma\in\mathit{I}_{r+}^{\boxplus}$ such that 
\begin{equation*}
\mathcal{C}_{\mu}^{\boxplus}(zS_{\mu}(z))=\mathcal{C}_{\sigma}^{\boxplus
}(z^{2}S_{\mu}^{2}(z)).
\end{equation*}
Since $z=\mathcal{C}_{\sigma}^{\boxplus}(zS_{\sigma}(z))$, we have that $%
S_{\sigma}(z)=zS_{\mu}^{2}(z)=zS_{\lambda}^{2}(z)S_{\mathrm{a}}^{2}(z).$
From (\ref{arcsinStrans}) $S_{\mathrm{a}}^{2}(z)=(z+2)/z$ and since $S_{%
\mathrm{m}_{2}}(z)=1/(z+2)$ we have $S_{\mathrm{m}_{2}}(z)S_{\sigma}(z)=S_{%
\lambda}^{2}(z).$ Then we conclude $\lambda\boxtimes\lambda=\mathrm{m}%
_{2}\boxtimes\sigma$.

Next assume that $\mu=\lambda\boxtimes\mathrm{a}$ and that there exist $%
\sigma\in\mathit{I}_{r+}^{\boxplus}$ such that $\lambda\boxtimes \lambda=%
\mathrm{m}_{2}\boxtimes\sigma$. From Theorem \ref{main1}, we can find $%
\tilde{\mu}\in\mathit{I}_{s}^{\boxplus}$ satisfying 
\begin{equation*}
\mathcal{C}_{\tilde{\mu}}^{\boxplus}(z)=\mathcal{C}_{\sigma}^{%
\boxplus}(z^{2}).
\end{equation*}
Then 
\begin{equation*}
\mathcal{C}_{\tilde{\mu}}^{\boxplus}(zS_{\tilde{\mu}}(z))=\mathcal{C}%
_{\sigma }^{\boxplus}(z^{2}S_{\tilde{\mu}}^{2}(z))=\mathcal{C}%
_{\sigma}^{\boxplus }(zS_{\sigma}(z))
\end{equation*}
from which we conclude that $zS_{\tilde{\mu}}^{2}(z)=S_{\sigma}(z)=(z+2)S_{%
\lambda}^{2}(z)$. Then $S_{\tilde{\mu}}^{2}(z)=S_{\mathrm{a}%
}^{2}(z)S_{\lambda}^{2}(z)$ and therefore $\tilde{\mu}=\lambda\boxtimes 
\mathrm{a}=\mu$.
\end{proof}

\begin{remark}
In the above theorem, if $\sigma\in\mathit{I}_{r+}^{\boxplus}$ is $\boxtimes 
$--2 divisible, $\lambda=\overline{\mathrm{m}_{2}}\boxtimes\overline{\sigma}$%
.
\end{remark}

\begin{example}
(1) If $\mathrm{b}_{s}$ is the symmetric beta $(1/2,3/2)$ distribution, $%
\mathrm{b}_{s}$ is the free type $AS$ distribution since $\mathrm{b}_{s}=%
\mathrm{m}\boxtimes\mathrm{a}$. We recall that $\mathrm{b}_{s}$ \ is the
free type $W$ distribution but not the free type $G$ distribution. Then the
class of free type $W$ distributions is a proper subclass of the class free
type $AS$. \newline
(2) The Wigner-free Poisson distribution (i.e. $\mu=\mathrm{w}\boxtimes%
\mathrm{m}$) is the free type $AS$ distribution. To see this observe that 
\begin{equation*}
S_{\mu}(z)=\frac{1}{\sqrt{z}(z+1)}=\sqrt{\frac{z+2}{z}}\frac{1}{z+1}\frac {1%
}{\sqrt{z+2}}=S_{\mathrm{a}}(z)S_{\mathrm{m}}(z)\frac{1}{\sqrt{z+2}},
\end{equation*}
where $\frac{1}{\sqrt{z+2}}$ is the S-transform of the probability measure $%
\overline{\mathrm{m}_{2}}$. Then $\mu=\mathrm{a}\boxtimes\mathrm{m}\boxtimes%
\overline{\mathrm{m}_{2}}$ and 
\begin{equation*}
(\mathrm{m}\boxtimes\overline{\mathrm{m}_{2}})\boxtimes(\mathrm{m}\boxtimes%
\overline{\mathrm{m}_{2}})=\mathrm{m}^{\boxtimes2}\boxtimes \mathrm{m}_{2}
\end{equation*}
and therefore $\mu$ is $\boxplus$--infinitely divisible and thus the free
type $AS$ distribution.
\end{example}

We finally show that the class type $AS$ distributions does not coincide
with the class \textit{I}$_{s}^{\boxplus}$ of all symmetric free infinitely
distributions.

\begin{example}
Let $\mu_{c}$ be as in example \ref{exNW}. Then 
\begin{equation}
S_{\mu_{c}}(z)=\sqrt{\frac{z+2}{z}}\frac{1}{\sqrt{(z+2)(z+2c)}}.
\label{aux3}
\end{equation}
If we take $c=1/15$, there is not a probability measure whose $S$--transform
is $\frac{1}{\sqrt{(z+2)(z+2c)}}$ by a similar argument as in Proposition %
\ref{No2divMP}. Since $\sqrt{(z+2)/z}$ is the $S$-transform of the arcsine
distribution $\mathrm{a}$ on $(-1,1),$ from (\ref{aux3}) and (\ref{DefMulCon}%
) we have that $S_{\mu_{c}}(z)$ cannot be the $S$-transform of a
multiplicative convolution with $\mathrm{a.}$
\end{example}

\section*{Acknowledgments}

The authors thank O. Arizmendi, M. Fevrier, W. M\l otkowski and J-C. Wang
for their stimulating remarks in the course of the preparation of this work.
They also thank the referees for their valuable and important comments that
have improved the presentation of the paper. 


\section*{Symbols}

\noindent $\mathcal{P}$ : the set of all probability measures on $\mathbb{R}$.%
\newline
$\mathcal{P}_{+}$ : the set of all probability measures on $\mathbb{R}_{+}$.%
\newline
$\mathcal{P}_{s}$ : the set of all symmetric probability measures on $\mathbb{%
R}$.\newline
$\mathcal{C}_{\mu }^{\ast }$ : the classical cumulant transform of
the probability measure $\mu $.\newline
$G_{\mu }$ : the Cauchy transform of the probability measure $\mu $.\newline
$F_{\mu }$ : the reciprocal of the Cauchy transform of the probability measure $%
\mu $.\newline
$\mathcal{C}_{\mu }^{\boxplus }$ : the free cumulant (or $R-$)transform of
the probability measure $\mu $.\newline
\textrm{$I$}$^{\ast }$ : the set of all (classical) infinitely divisible
distributions on $\mathbb{R}$.\newline
\textrm{$I$}$_{+}^{\ast }$ : the set of all positive (classical) infinitely
divisible distributions on $\mathbb{R}_{+}$.\newline
$I_{s}^{\ast }$ : the set of all symmetric (classical) infinitely divisible
distributions on $\mathbb{R}$.\newline
\textrm{$I$}$^{\boxplus }$ : the set of all free infinitely divisible
distributions on $\mathbb{R}$.\newline
$I_{s}^{\boxplus }$ : the set of all free symmetric infinitely divisible
distributions on $\mathbb{R}$.\newline
$I_{r+}^{\boxplus }$ : the set of all free regular infinitely divisible
distributions on $\mathbb{R}_{+}$. \newline
$\mathcal{L}(X)$ : the law of a real random variable $X$.\newline
$\nu _{\mu }$ : the L\'{e}vy measure of $\mu \in $\textrm{$I$}$^{\ast }$ or 
\textrm{$I$}$^{\boxplus }$.\newline
$S_{\mu }$ : the S-transform of the probability measure $\mu \in \mathcal{P_{+}}$%
. \newline
$\mathrm{w}_{b,a}$ : the Wigner (or semicircle) distribution with mean $b$
and variance $a$. \newline
$\mathrm{w}$ : the Wigner (or semicircle) distribution with mean $0$ and
variance $1$. \newline
$\mathrm{m}_{c}$ : the Marchenko-Pastur (or free Poisson) distribution with
parameter $c>0$. \newline
$\mathrm{m}$ : the Marchenko-Pastur (or free Poisson) distribution with
parameter $1$. \newline
$\mathrm{a}_{s}$ : the symmetric arcsine distribution on $(-s,s)$. \newline
$\mathrm{a}$ : the symmetric arcsine distribution on $(-1,1)$. \newline
$\mathrm{a}_{s}^{+}$ : the positive arcsine distribution on $(0,s)$. \newline
$\mathrm{a}^{+}$ : the positive arcsine distribution on $(0,1)$. \newline
$\gamma _{b,a}$ : the Gaussian distribution with mean $b$ and variance $a$. 
\newline
$\gamma _{s}$ : the Gaussian distribution with mean $0$ and variance $s$. 
\newline
\textrm{p}$_{c}:$ the classical Poisson distribution with mean $c>0.$ 
\newline
$\gamma _{s}^{(2)}$ : the gamma distribution with shape parameter $1/2$ and
scale parameter $s$.

\end{document}